\documentclass[11pt]{imsart}
\usepackage{amssymb}
\pdfoutput=1

\usepackage[margin=3cm]{geometry}

\usepackage{amsmath, moreverb}
\usepackage{makecell}
\usepackage{color}
\usepackage[round]{natbib}
\usepackage{amsfonts, fancybox}
\usepackage{amssymb}
\usepackage[american]{babel}
\usepackage{graphicx,pstricks, pst-node, color}
\usepackage{psfrag}\usepackage{subfigure}
\usepackage{flafter, bm, dsfont}
\usepackage[section]{placeins}

\startlocaldefs

\newcommand{\lnormpi}[1]{\ensuremath \Vert #1 \Vert_\pi}
\newcommand{\ippi}[2]{\ensuremath \langle #1, #2 \rangle_\pi}
\newcommand{\modulo}{\ensuremath \mathrm{mod}~}
\newcommand{\eye}{\mathbf{I}}
\newcommand{\reals}{\mathbf{R}}

\newcommand{\cN}{\mathcal{N}}

\newcommand{\cX}{\mathcal{X}}

\newcommand{\KL}{{\sf KL}}

\newcommand{\R}{\mathbf{R}}
\newcommand{\fP}{\mathbf{P}}

\newcommand{\E}{\mathbf{E}}

\makeatletter
\def\namedlabel#1#2{\begingroup
   \fPtect\def\@currentlabel{#2}%
   \label{#1}\endgroup
}
\makeatother

\usepackage{amsthm}
\usepackage{hyperref}
\usepackage{cleveref}
\hypersetup{colorlinks=true,citecolor=blue,urlcolor=blue}
\newcommand{\BlackBox}{\rule{1.5ex}{1.5ex}}
\renewenvironment{proof}{\par\noindent{\bfseries\upshape
  Proof\ }}{\hfill\BlackBox\\[2mm]}

\newtheorem{theorem}{Theorem}

\newtheorem{lemma}[theorem]{Lemma}
\newtheorem{proposition}[theorem]{Proposition}
\newtheorem{remark}[theorem]{Remark}

\newtheorem{definition}[theorem]{Definition}

\makeatletter
\newcommand*{\defeq}{\mathrel{\rlap{%
                     \raisebox{0.3ex}{$\m@th\cdot$}}%
                     \raisebox{-0.3ex}{$\m@th\cdot$}}%
                    =}
\newcommand*{\eqdef}{=
  \mathrel{\rlap{%
      \raisebox{0.3ex}{$\m@th\cdot$}}%
    \raisebox{-0.3ex}{$\m@th\cdot$}}%
}
\makeatother

\crefname{theorem}{Theorem}{Theorems}
\crefname{observation}{Observation}{Observations}
\crefname{claim}{Claim}{Claims}
\crefname{condition}{Condition}{Conditions}
\crefname{example}{Example}{Examples}
\crefname{fact}{Fact}{Facts}
\crefname{lemma}{Lemma}{Lemmas}
\crefname{corollary}{Corollary}{Corollaries}
\crefname{definition}{Definition}{Definitions}
\crefname{remark}{Remark}{Remarks}

\newcommand{\makemark}[1]{#1}

\runauthor{\sc Berthet and Kanade}
\runtitle{Statistical Window Testing for the Initial Distribution}
\begin{document}

\begin{center}
{\LARGE Statistical Windows in Testing for the Initial Distribution\\ of a Reversible Markov Chain} 
\vskip 0.4cm
{Quentin Berthet \footnote{Statistical Laboratory, DPMMS, University of Cambridge} and Varun Kanade \footnote{Department of Computer Science, University of Oxford}}
\vskip 1.5 cm
	\small{\textbf{Abstract}}
\end{center}

\maketitle

\begin{small} 
	We study the problem of hypothesis testing between two discrete distributions, where we only have access to samples after the action of a known reversible Markov chain, playing the role of noise. We derive instance-dependent minimax rates for the sample complexity of this problem, and show how its dependence in time is related to the spectral properties of the Markov chain. We show that there exists a wide {\em statistical window}, in terms of sample complexity for hypothesis testing between different pairs of initial distributions. We illustrate these results in several concrete examples.
\end{small}

\thispagestyle{empty}

\section{Introduction}
\label{sec:intro}

Random walks on graphs, or Markov chains more generally, have long served as a
natural model for data observed through a noisy channel. In many applications, one wishes to determine the origin of a random walk, after a certain number of steps: Where has a rumor started in a social network? What is the distribution of a deck of cards before shuffling? What are the ancestors of current species before evolution? What is the initial configuration of spins in Glauber dynamics? In these cases, the initial distribution is considered as the true information, and observations are made after the action of a few steps of the Markov
chain. In this work, we consider the problem of testing for the initial distribution: determining which of two candidate distributions $\mu$ and $\mu'$ is the initial distribution, based on an i.i.d. sample after $t$ steps. The mixing time of the chain can
be interpreted as the time at which ``all'' starting information is lost.  
However, depending on the two candidate distributions over the starting states, this total
loss of information may occur at a time much sooner than the mixing time. In
this work, we characterize precisely this rate of loss of information in terms of spectral
properties of the transition matrix: we give a theoretical guarantee on the performance of a test in terms of necessary sample size as a function of all the parameters of the problem, and show that this dependency is tight, using information-theoretic tools. In particular, we show that the sample complexity of the hypothesis testing problem between $\mu$ and $\mu'$, and its dependency in time $t$, depends critically on the pair $(\mu,\mu')$ in an explicit manner. We call this wide range of sample complexities the {\em statistical window}. Pairs of distributions that exhibit behaviour at the extreme ends of this statistical window can be explicitly constructed using the spectrum of the associated Markov chain and we illustrate this phenomenon on several concrete examples.

Recovering information about a discrete distribution with access to samples is one of the central problems of statistical theory, going back at least to \cite{Lap12}. This essential problem has attracted much attention in the modern treatment of learning theory, on problems related to learning, testing and estimation. Recently, there has been renewed focus on learning discrete distributions from a sample---typically, it is assumed that the (unknown) distribution satisfies certain properties such as $k$-modality or monotonicity, which, in certain cases, allows significantly improved sample complexity over the basic approach of using the empirical distribution (see e.g. \cite{chan2013learning,daskalakis2012learning,pmlr-v40-Kamath15,DiaKanNik14,daskalakis2015learning,diakonikolas2016learning,DiaKan16,DiaKanNik17,DiaGouPee17,ValVal14,valiant2016instance}).  A related area of research is that of property testing, i.e. to test whether a distribution satisfies some property such as uniformity or monotonicity, from a sample (see e.g.~\cite{CDVV:2014,canonne2014testing,diakonikolas2015optimal,valiant2011testing}).  Yet another area of interest has been estimating quantities related to a discrete distribution, such as support size or entropy (see e.g.~\cite{VV:2011,acharya2014complexity,WuYan14,WuYan16,OrlSurWu16}).

In much of this literature, it is often assumed that one has direct access to independent samples from the true unknown distribution of interest $\mu$.%
\footnote{There have been some recent advances where some of these problems can be solved even when some (small) fraction of the data has been adversarially tampered with.}

As stated above, we consider in this work a setting where we only have access to these samples after a reversible Markov chain has acted on them $t\ge0$ times. In doing so, we allow for the introduction of noise, in the form of the action of a known Markov chain with transition matrix $P$. Formally, this is equivalent to learning about $\mu$, with access to $\mu_t = \mu P^t$, and can be seen as a statistical inverse problem. In our setting, as $t$ increases and $\mu_t$ approaches the stationary distribution, more information is lost and the statistical problem becomes more difficult. This is in stark contrast to the usual applications of Markov chains in statistical learning, in particular for Markov Chain Monte Carlo methods, where the stationary distribution $\pi$ is the quantity of interest, from which it is difficult to sample, and a large $t$ is desirable. In our setting, $t$ can be understood as a way to measure the amount of noise, and we seek to understand how the difficulty increases with it. This is a common point of view in some continuous settings: if $\mu = \delta_x$ for some $x \in \R$, the action of the heat kernel for time $t$ leads to a distribution $\mu_t = \cN(x,t)$, and the impact of $\sigma^2=t$ on the statistical difficulty of recovering $x$ is clear. We transfer this idea to discrete distributions, and the action of a Markov chain is the most natural way to introduce noise.

\begin{figure}[h!]
\label{FIG:histog_markov}
    \begin{center}
      \includegraphics[width= 0.8\textwidth]{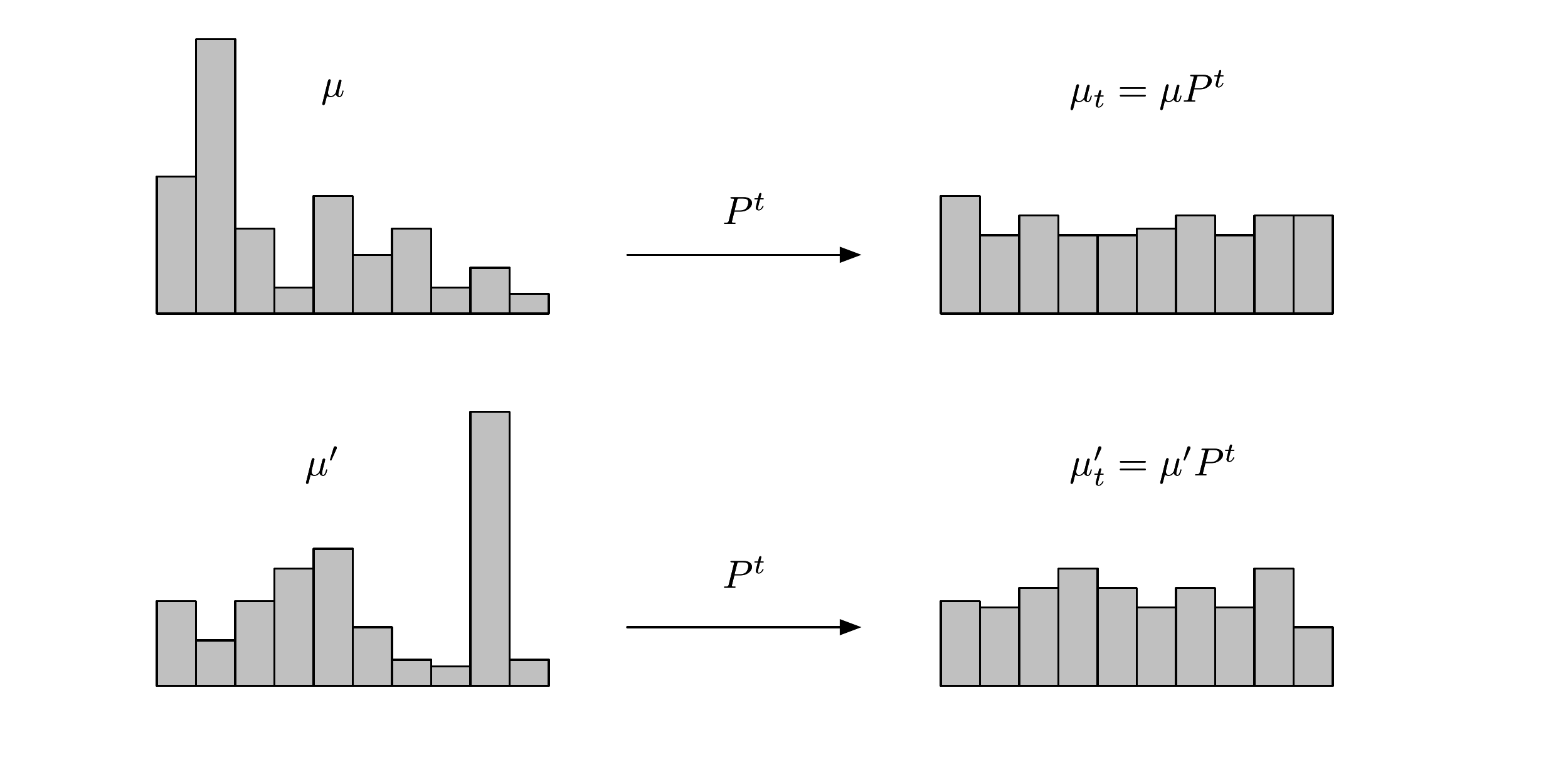}
			\caption{The hypothesis testing problem between $\mu$ and $\mu'$ is easier if we have direct sample access to these distributions. It is harder with the action of $P^t$, as the distributions look increasingly alike. For very large $t$, the distribution is very close to stationarity and all initial information is lost.}
    \end{center}
\end{figure}

Formally, we focus here on the fundamental problem of hypothesis testing between
two known distributions $\mu$ and $\mu'$, based on samples from $\mu_t = \mu
P^t$ or $\mu'_t = \mu' P^t$. This choice allows us to illustrate, for many
natural examples, the impact of the pair $\mu,\mu'$ on the instance-dependent sample
complexity of this problem, the required sample size $n^*_{\mu,\mu',t}$ to
solve this problem with a small probability of error, and in particular its dependency in time. There is much focus in
the field of mathematical statistics on the analysis of the maximum probability
of error in hypothesis testing. It is well-understood that it is equivalent to
studying the {\em total variation distance} $d_{\sf TV}(\mu_t^{\otimes
n},\mu_t^{\prime \otimes n})$ between the two distributions of samples of size
$n$. This quantity grows in $n$ and decays in $t$, and the goal of our analysis
is to describe the trade-off between these two phenomena: we establish how
large $n$ and how small $t$ need to be for the total variation distance to be
bounded away from 0, for testing to be possible. We are particularly interested in showing that the behaviour in $t$ is far from universal: we exhibit pairs $(\mu,\mu')$ for which it is very different, as well as a systematic manner to construct them.

Our analysis bridges two fields where studying this quantity is a central
problem. On the one hand, for fixed $t$, given $\mu_t$ and $\mu'_t$,
understanding the growth in $n$ of this distance is one of the central problems
of mathematical statistics and information theory. It is difficult to
establish the growth of $d_{\mathsf{TV}}(\mu_t^{\otimes n}, \mu_t^{\otimes n})$
in terms of $n$ and $d_{\mathsf{TV}}(\mu_t, \mu_t)$ directly.

In our analysis, this quantity is controlled by comparison to other notions of ``distances''.
On the other hand, one of the most important questions in the study of Markov chains is the behavior of $\mu_t$, and in particular the decay of its total variation distance to the stationary distribution, as a function of $t$. One key quantity is that of mixing time, which describes the time $t$ at which the total variation goes to $e^{-1}$.

After mixing, the total variation distance to the stationary distribution
$\pi$ decays exponentially. The mixing time represents the time at which the
sample complexity of the hypothesis testing problem explodes, {\em for
all} pairs of initial distributions $\mu,\mu'$.  What emerges from our analysis
in Section~\ref{sec:samplecomp} is that for $t$ less than the mixing time, this is
far from universal over all pairs $\mu,\mu'$, and that in many natural Markov
chains, the sample complexity can vary dramatically. In particular, asking
when $\mu_t$ and $\mu'_t$ are close is a very different question from asking
when $\mu_t$ is close to the stationary distribution $\pi$. In particular, we construct explicitly examples of pairs of initial distributions $(\mu,\mu')$ whose sample complexity has a very different behavior in time.

Our main results are in Section~\ref{sec:sample}, where we show that the instance sample complexity is
\[
n^*_{\mu,\mu',t} \asymp 1/\sum_{i=2}^d \lambda_i^{2t} (\alpha_i -\alpha'_i)^2\, ,
\]
where the $\lambda_i$ are the eigenvalues of $P$ and the $\alpha_i,\alpha'_i$ are components of $\mu$ and $\mu'$ along its eigenvectors. 
In particular, a very rich structure emerges, where all the eigenvalues play a role: depending on the values of $(\alpha_i-\alpha'_i)$, the sample complexity as a function of time can be driven by terms of order $\lambda_i^{2t}$ for any $i$, and does not only depend on the spectral gap. At the extremes, we describe the {\em statistical window}, the ratio of sample complexities $n^*_{\mu,\mu',t}/n^*_{\gamma, \gamma',t}$, for pairs $\mu,\mu'$ and
$\gamma,\gamma'$ of initial distributions with comparable initial total sample complexity. We show that 
	\[
		\frac{n^*_{\mu,\mu',t}}{n^*_{\gamma,\gamma',t}} \asymp
		\frac{n^*_{\mu,\mu',0}}{n^*_{\gamma,\gamma',0}} \left(
		\frac{\lambda_{[2]}}{\lambda_{[d]}}\right)^{2t} \, ,
	\]

is governed by the ratio between the eigenvalues of the Markov chain with the largest and smallest absolute value (see Theorem~\ref{THM:window}). We illustrates these findings by deriving the sample complexity for
several concrete examples of Markov chains in Section~\ref{sec:applications}.
All original proofs are given in the appendix.

\paragraph{Related Work.}  There has been considerable work regarding
reconstruction of a signal observed through noisy channels. An area of
particular interest has been information flow on (rooted) trees, where the
label (color) of the root of a (possibly infinite) tree is chosen according to
a discrete distribution.  Each edge of the tree acts a noisy channel, given by
the transition matrix of a Markov chain $P$. The goal is to reconstruct the
signal at the root, given the information at the leaves. Sharp results are
known for several cases depending on the branching factor of the tree and the
eigenvalues of $P$~\citep{evans2000broadcasting,mossel2003information,mossel2001reconstruction}.
Although, the noise model is similar to the one we consider in
this paper, the goal is significantly different---their focus is on
reconstruction of a single signal from several (possibly correlated) corrupted
observations. Other problems based on similar principles include
population recovery \citep{DviRaoWig12,PolSurWu17} and trace reconstruction
\citep{Lev01,McGPriVor14,HarHolPer17,HolPemPer18}.

A recent line of work has focused on testing whether a Markov chain $P$ is
identical to a fixed chain $P^\prime$ or sufficiently far from it, given a
single trajectory $X_0, \ldots, X_t$ generated from $P$~\citep{DDG:2017}. Their
work does not assume the knowledge of $P$. In contrast, for the testing problem
considered in this work, it is easy to see that in the absence of some
knowledge of $P$, the testing problem is impossible, e.g. one can easily
construct pairs of starting distribution and transition matrix, $(\mu, P)$ and
$(\mu^\prime, P^\prime)$, such that the distributions $\mu P$ and $\mu^\prime
P^\prime$ are identical. Also, there is little to be gained by observing the
trajectory in our setting as all the relevant information is contained in the
first observation of the trajectory. Other statistical problems based on the observations from a Markov Chain include learning a graphical model from Glauber dynamics \citep{BreGamSha14}, the mixing time \citep{HsuKonSze15}, or the entropy rate \citep{KamVer16}. The notions of statistical distances in relation with Markov chains, in particular with respect to their stationary distributions have also been explored in \citep{BreNag17}. The problem of how initial information is lost with the action of a Markov chain is also considered in \citep{GolBrePol18}, who study the case of Glauber dynamics in the context of information storage.

\paragraph{Notation.} Throughout the paper, $\delta$ is the probability of
error and $\varepsilon \in (0, 1)$ is a measure of \emph{non-degeneracy} of
distributions. The notation $\asymp$ indicates equality up to factors that may
depend only on $\delta$ and $\varepsilon$. Standard notions from
information theory which we use are defined in Appendix~\ref{APP:notation}.

\section{Problem description}
\label{sec:problem_desc}
Let $\cX$ be a set of size $d$, and $P$ the transition matrix of a known
irreducible reversible Markov chain on $\cX$. We observe $X_1,\ldots,X_n$ that
are i.i.d. draws from this Markov chain after $t \ge 0$ steps, with an unknown
{\em initial distribution} $\nu$. The distribution of the $X_i$ is therefore
$\nu_t = \nu P^t$. Our objective is to determine, given two distributions $\mu$
or $\mu'$ on $\cX$, which one of these is the initial distribution $\nu$, based
on the observation of a sample $(X_i)_{i \in [n]}$. This is equivalent to a hypothesis
testing problem between $\mu_t = \mu P^t $ and $\mu'_t =  \mu' P^t$, for known
$t$ and $P$, based on an i.i.d. sample of size $n$.

For any test $\psi : \cX^n \to \{\mu,\mu'\}$, its performance is measured in terms of its {\em maximum probability of error}
\[
\max_{\nu_0 \in \{\mu,\mu'\}} \fP^{\otimes n}_{\nu_t}(\psi \neq \nu_0) = \fP^{\otimes n}_{\mu_t}(\psi \neq \mu) \vee \fP^{\otimes n}_{\mu'_t}(\psi \neq \mu')\, .
\]
In this work, we analyze the {\em sample complexity} of this problem, i.e.
the required sample size $n^*_{\mu, \mu',t}$ to have a small probability of
error. Formally, for $n \gtrsim n^*_{\mu, \mu',t}$, we have that the
probability of error is smaller than $\delta$ for some test, and for $n \lesssim
n^*_{\mu, \mu',t}$, that it is greater than $1/2-\delta$ for all tests, for a
fixed probability $\delta \in (0,1/4)$, up to multiplicative constants. If $\mu_t=\mu'_t$ for some $t$, the statistical problem is impossible and we write $n^*_{\mu, \mu',t} = \infty$.

As discussed in the introduction, the maximum probability of error is related
to the total variation distance $d_{\sf TV}(\mu_t^{\otimes n},\mu_t^{\prime \otimes n})$, and our analysis
relies on understanding the behavior of this quantity. We analyze its growth in
$n$ using tools of mathematical statistics, and its decay in $t$ through the lens of
spectral analysis for reversible Markov chains in Section~\ref{sec:sample}.

We use the following notion of two distributions having a bounded-likelihood ratio and state an immediate consequence below.

\begin{definition}[Bounded likelihood-ratio]
	\label{AS:initial}
	Two distributions $\mu$ and $\mu'$ have an $\varepsilon$-bounded likelihood-ratio if for all $x\in \cX$, $\varepsilon\le \mu_x/\mu'_x
	\le 1/\varepsilon$.
\end{definition}

\begin{proposition}
\label{PRO:inherit}
If $\mu,\mu'$ have an $\varepsilon$-bounded likelihood-ratio, so do $\mu_t=\mu P^t$ and $\mu'_t = \mu' P^t$
\end{proposition}

We make the assumption in \makemark{some of our results}, that the three
distributions, $\pi$, the stationary distribution of the Markov chain $P$, and
$\mu$, $\mu^\prime$, the initial distributions, pairwise have an
$\varepsilon$-bounded likelihood-ratio; this is to avoid some pathological cases where
the statistical complexity is very small, e.g. if $\mu_t$ and $\mu'_t$ do not
have the same support, some observations suffice to solve the testing problem
with probability of error $0$. This assumption therefore ensures that it cannot
be a trivial problem. For the Markov chains that we consider, this
property is eventually satisfied for $t$ larger than some fixed quantity,
provided $P$ is aperiodic. 

Our main interest is to exhibit the statistical window phenomenon: we exhibit
pairs of initial distributions for which the hypothesis testing problem has
very different sample complexities. These distributions all satisfy the
assumption above. Further, we show in Section~\ref{SEC:general} that these results can
be extended to more general cases: if the distributions $\mu, \mu^\prime$ do
not satisfy the \makemark{bounded likelihood-ratio assumption
(Definition~\ref{AS:initial})}, the hypothesis testing question can be
rephrased in terms of distributions that do, up to a loss in multiplicative
factors.

\section{Sample complexity}
\label{sec:samplecomp}

\subsection{Reversible Markov chains}

\label{sec:markov}

In this problem, we have access to a sample of size $n$ from either $\mu_t$ of $\mu'_t$. It is therefore important to understand the behavior of these two distributions as $t$ increases, and in particular how quickly they become similar, depending on the initial starting points $\mu$ and $\mu'$. We state some basic properties of reversible Markov chains and their spectra
without proof; these can be found in standard texts on Markov chains (e.g.
\cite{LPW:2008}). Recall that $P$ is the transition matrix of an irreducible
reversible Markov Chain on a finite state space $\cX$ with stationary
distribution $\pi$. We adopt the convention that $P_{ij}$ is the probability of
transitioning from state $i$ to $j$; as $P$ is reversible we have $\pi_i P_{ij}
= \pi_j P_{ji}$ for all $i, j$. Denote by $\Pi$ the diagonal matrix with
$\Pi_{ii} = \pi_i$ and consider the matrix $Q := \Pi^{\frac{1}{2}} P \Pi^{-
\frac{1}{2}}$. As $P$ is reversible and $\pi$ is the stationary distribution,
$Q_{ij} =  \sqrt{\frac{\pi_i}{\pi_j}} P_{ij} = \sqrt{\frac{\pi_j}{\pi_i}}
P_{ji} = Q_{ji}$, and as a result $Q$ is symmetric.  The following proposition
holds for reversible Markov Chains.

\begin{proposition}
	Let $\nu_1, \ldots, \nu_d$ be the eigenvectors\footnote{As $Q$ is symmetric the left and right eigenvectors are the same.} and $1 = \lambda_1 >
	\lambda_2 \geq \lambda_3 \geq \cdots \geq \lambda_d \geq -1$ the corresponding eigenvalues
	of $Q$. Then,\\
	\noindent
		$-$ $v_i = \Pi^{-\frac{1}{2}} \nu_i$ is a right eigenvector of $P$ with eigenvalue $\lambda_i$; for $\lambda_1 = 1$, $v_1 = (1, \cdots, 1)^\top$. \\
		\noindent
		$-$ $u_i = \Pi^{\frac{1}{2}} \nu_i = \Pi v_i$ is a left eigenvector of $P$ with eigenvalue $\lambda_i$; for $\lambda_1 = 1$, $u_1 = \pi$.
\end{proposition}

We consider $1 = \lambda_{[1]}, \lambda_{[2]}, \ldots, \lambda_{[d]}$ as the
ordering of the eigenvalues by their absolute values, i.e. $1 = |\lambda_{[1]}|
\geq |\lambda_{[2]}| \geq \cdots \geq |\lambda_{[d]}| \geq 0$.  It is worth
pointing out that $\lambda_{[2]} = \lambda_d = -1$ is possible if and only if
$P$ is periodic with period $2$, i.e. the underlying graph on $\cX$ is
bipartite. It is common to consider \emph{lazy} chains, where one considers
the transition matrix $(1 - q) P + q \eye$. An eigenvalue $\lambda$ of the
original chain yields an eigenvalue $(1 - q) \lambda + q$ for the lazy chain.
In particular, for $q \geq 1/2$, all eigenvalues become non-negative. We
discuss the impact of adding \emph{laziness} to the problem of testing initial
distributions in Section~\ref{sec:sample}.

\paragraph{Inner Product and Norm with respect to $\pi$.} As $P$ is irreducible, $\pi_x > 0$ for every $x \in \cX$. Thus, we can define an inner product over $\reals^d$ as follows:
\begin{align}
	\langle u, u^\prime \rangle_{\pi} = \sum_{x} \frac{u_x u_x^\prime}{\pi_x} \label{eqn:pi-inner-prod}
\end{align}

We denote the associated norm as $\lnormpi{u} = \sqrt{\ippi{u}{u}}$. The
following lemma states that the left eigenvectors form an orthonormal basis
with respect to the inner product $\ippi{\cdot}{\cdot}$. 

\begin{lemma} 
	\label{lem:markov-basis}
	The left eigenvectors $u_1, \ldots, u_d$ of $P$ form an orthonormal basis
	with respect to the inner product $\ippi{\cdot}{\cdot}$. Furthermore, for
	any $u \in \reals^d$ with $\sum_x u_x = 1$, we have $\ippi{u}{\pi} = 1$; in
	particular for $u_1 = \pi$, we have $\lnormpi{\pi} = 1$. By orthogonality,
	we also have that for $i \geq 2$, $\sum_{x} u_{i, x} = \ippi{u_i}{\pi} = 0$;
	henceforth without loss of generality, we will assume that they form an
	\emph{orthonormal} basis, i.e. $\lnormpi{u_i} = 1$ for each $i$.
\end{lemma}

The distance $\lnormpi{\mu - \mu^\prime}$ might seem strange at first sight.
Note that $\lnormpi{\mu - \mu^\prime}^2 = \sum_{x} \frac{(\mu_x -
\mu_x^\prime)^2}{\pi_x}$. Let us consider the case when one of the two
distributions, say $\mu^\prime$ is the stationary distribution $\pi$: then, it is simply the $\chi^2$ divergence, $D_{\chi^2}(\mu, \pi)$ between $\mu$
and $\pi$. In our analysis, we compare this distance to other notions of distances between $\mu_t$ and $\mu'_t$, to obtain guarantees on the sample complexity of this problem. We rely on the spectral properties of the transition matrix to understand the temporal evolution of an initial distribution over $\cX$. 
\subsection{Guarantees for likelihood-ratio test}

\label{sec:infotheory_desc}

For this testing problem, we show guarantees for the performance of the likelihood ratio test. It is based on the  log-likelihood ratio statistic $L_n$ between $\mu_t$ and $\mu'_t$, given by
\[
L_n = \sum_{x \in \cX} \hat \mu_{t,x} \log(\mu_{t,x}/\mu'_{t,x}) = \langle \ell, \hat \mu_t \rangle\, ,
\]
where $\hat \mu_t$ is the empirical or observed distribution of the $(X_i)_{i \in [n]}$, and $\ell_x = \log(\mu_{t,x}/\mu'_{t,x})$.
\begin{definition}
The likelihood ratio test $\psi_{\sf LR}$ takes $(X_i)_{i \in [n]} \in \cX^n$ as input and outputs $\mu$ or $\mu'$ such that
\[
\psi_{\sf LR} = \left \{
\begin{array}{rl}
\mu & \text{if } L_n>0\, ,\\
\mu' &\text{if } L_n\le 0\, .
\end{array} \right.
\]
\end{definition}

The Kullback-Leibler divergences between $\mu_t$ and $\mu'_t$ are naturally associated to the quantity $L_n$, as its expected value under these distributions. This divergence is a statistical measure of divergence that captures well the sample complexity of the problem, and also appears in large deviations, in the description of the asymptotic behavior of $\hat \mu_t$. The connections between notions of ``distances'' between distributions and sample complexity have been extensively studied, see e.g. \citep{PolWu17} and references therein, and \citep{BerHarKon14}.

\subsection{Sample Complexity Guarantees}

\label{sec:sample}

\begin{theorem}
\label{THM:main1}
For two initial distributions $\mu,\mu'$, with $\mu, \mu^\prime, \pi$ all pairwise having $\varepsilon$-bounded likelihood-ratios for some $\varepsilon\in(0,1)$, the likelihood-ratio test $\psi_{\sf LR}$ has probability of error less than $\delta$ if 
\[
n \ge  C(\varepsilon,\delta)/\sum_{i=2}^d \lambda_i^{2t} (\langle u_i,\mu \rangle_\pi - \langle u_i,\mu' \rangle_\pi)^2\, ,
\]
for $C(\varepsilon,\delta) = 16 \varepsilon^{-5/2} \log(1/\delta)$.
\end{theorem}
This result is not obtained by focusing on the behavior of the random variable $L_n$ and using concentration inequalities, as is usually the case in such problems, but directly by analyzing and linking several notion of distances between the distributions $\mu^{\otimes n}_t$ and $\mu^{\prime \, \otimes n}_t$. Indeed, this allows us to understand simultaneously the growth in $n$ and convergence in $t$ of these distances and to give guarantees on how large $n$ needs to be for any fixed $t$. \makemark{We present the alternate point of view, more common in the analysis of Markov chains, in Section~\ref{SEC:time} below.} Furthermore, using a different analysis of other measures of statistical distance between distributions, we show that this guarantee on the performance of the likelihood ratio test is optimal: up to constants, it is tight for {\em all} tests depending only on the observation of a sample of size $n$. 
\begin{theorem}
\label{THM:main2}
For two initial distributions $\mu,\mu'$, with $\mu, \mu^\prime, \pi$ all pairwise having $\varepsilon$-bounded likelihood-ratios for some $\varepsilon\in(0,1)$, all tests have probability of error greater or equal to $1/2-\delta$ if
\[
n \le  c(\varepsilon,\delta)/\sum_{i=2}^d \lambda_i^{2t} (\langle u_i,\mu \rangle_\pi-\langle u_i,\mu' \rangle_\pi)^2\, ,
\]
for $c(\varepsilon,\delta) = 8\varepsilon \delta^2$.
\end{theorem}
When the sample size is smaller than this bound, no test can accurately identify the correct distribution, and significantly outperform a coin flip. Together, these results give a complete picture of the statistical complexity of
the hypothesis testing problem defined in Section~\ref{sec:problem_desc}. The sample complexity of this problem is of order
\[
n^*_{\mu,\mu',t} \asymp 1/\sum_{i=2}^d \lambda_i^{2t} (\langle u_i,\mu \rangle_\pi-\langle u_i,\mu' \rangle_\pi)^2\, .
\]
This expression gives a very clear understanding of how the initial
information is lost over time. The component of the difference between $\mu$ and $\mu'$ (seen as vectors in $\R^\cX$) 
aligned with eigenvectors with eigenvalues close to $0$ will be lost fast, while
that along those with eigenvalues close to $-1$ and $1$ will be
retained longer. As a consequence, different pairs of initial distributions have
very different statistical complexities. The results above allow us to describe exactly this phenomenon. The range of sample complexities for this problem can therefore be very large, and is governed by the spectral properties of the matrix.
\begin{theorem}
\label{THM:window}
	There are pairs of initial distributions $\mu,\mu'$ and $\gamma,
	\gamma'$, with $\mu, \mu^\prime, \pi$ and $\gamma, \gamma^\prime, \pi$
	pairwise having $\varepsilon$-bounded likelihood-ratios for some $\varepsilon\in(0,1)$,
	such that
	\[
		\frac{n^*_{\mu,\mu',t}}{n^*_{\gamma,\gamma',t}}  \asymp
		\frac{n^*_{\mu,\mu',0}}{n^*_{\gamma,\gamma',0}} \left(
		\frac{\lambda_{[2]}}{\lambda_{[d]}}\right)^{2t} \, .
	\]
\end{theorem}

In particular, if the initial statistical complexities $n^*_{\mu,\mu',0}$ and
$n^*_{\gamma,\gamma',0}$ are similar, the ratio scales like
$(\lambda_{[2]}/\lambda_{[d]})^{2t}$, and we refer to it as the {\em
statistical window}. If there is an eigenvalue that is negative, we can
arbitrarily increase the size of the statistical window with laziness, moving the
eigenvalue closer to $0$, or even make it infinite.
We analyze this window for examples of Markov chains in
the following section, describing as well the extremal pairs of initial
distributions at the two ends of this statistical window: the hypothesis
testing problems that become hard quickly, and those which are the least
affected by the action of $P^t$. In many applications, this provides an intuitive understanding of the type of questions that become hard, for several natural random processes, through the lens of loss of information.

\subsection{Guarantees without likelihood-ratio bounds}
\label{SEC:general}

Our main message is that there {\em exist} pairs of distributions whose associated hypothesis testing problems have vastly different sample complexity, and in particular that this phenomenon can be exhibited by taking distributions satisfying a bounded likelihood-ratio assumption. To analyze the sample complexity for two distributions that do not satisfy this assumption, one can reduce the problem to the case of two distributions that do.
\begin{definition}
For any $\eta \in (0,1)$, $\mu,\mu'$ distributions on $\cX$ and a Markov chain $P$ with stationary distribution $\pi$, we consider 
\[
	\beta=(\mu+\mu'+\pi)/3\, ,
\]
the average of $\mu,\mu'$, and $\pi$. The centered versions $\tilde \mu$ and $\tilde \mu'$ are defined for any $\eta \in (0,1)$ as
\[
	\tilde \mu = (1 - \eta)\mu + \eta \beta\, , \quad \tilde \mu' = (1 - \eta) \mu^\prime + \eta \beta\, .
\]
\end{definition}

Considering the hypothesis testing problem between $\tilde \mu$ and $\tilde \mu'$ only makes the statistical problem harder: it can be interpreted as drawing each sample point from $\beta$ instead of either $\mu$ or $\mu'$, with probability $\eta$. Note that $\tilde \mu, \tilde \mu'$, and $\pi$ all pairwise having $\eta/3$-bounded likelihood-ratios. Using these distributions, we generalize Theorem~\ref{THM:main1} as follows.

\begin{theorem}
\label{THM:general}
For two initial distributions $\mu,\mu'$, the likelihood-ratio test $\psi_{\sf LR}$ has probability of error less than $\delta$ if 
\[
n \ge  c \log(1/\delta)/\sum_{i=2}^d \lambda_i^{2t} (\langle u_i,\mu \rangle_\pi - \langle u_i,\mu' \rangle_\pi)^2\, ,
\]
for a universal constant $c>0$.
\end{theorem}
This result rests on the proof of Theorem~\ref{THM:main1} describing the sample complexity for the pair $(\tilde \mu,\tilde \mu')$, and controlling the difference in sample complexity with testing for the pair $(\mu,\mu')$.

However, the sample complexity could be much smaller: if $\mu$ and $\mu'$ have different supports, if this property still holds for $\mu_t$ and $\mu'_t$, the sample complexity can be of the order of a constant. Outside of such degenerate cases, if they have full supports, the lower bound of Theorem~\ref{THM:main2} can be recovered up to a multiplicative factor of the maximum of $\pi_x/\mu'_{t,x}$ - which is always finite, by following the same proof. Both upper and lower bounds are therefore valid up to constants if one of the two distributions is $\pi$, and the other has the same support. This last case allows also to showcase the full width of the statistical window, by taking distributions $\mu$ such that $\mu-\pi$ is aligned along different left eigenvectors of $P$.

\subsection{Statistical time guarantees}

\label{SEC:time}
Our results are presented in a fixed time, fixed probability setting, and we give guarantees in terms of how large the sample size $n$ needs to be. However, in most of the literature on Markov chains, there is no sample size, and results are given in terms of guarantees on the time $t$. Some of our results can be formulated in a similar manner: given a fixed sample size $n$, what is the {\em statistical time} $t^*_{\mu,\mu',n}$ such that the testing problem is possible when $t  \le t^*_{\mu,\mu',n}$ and becomes impossible when $t \ge t^*_{\mu,\mu',n}$, up to terms involving only $\delta$ and $\varepsilon$. There is no general expression for this statistical time, however it can be made explicit for many pairs of initial distributions $\mu,\mu'$, also a direct consequence of Theorem~\ref{THM:main1} and~\ref{THM:main2}.

\begin{theorem}
\label{THM:time}
For any reversible Markov chain $P$, and for every $i \in [d]$, there exists a pair of initial distributions $\mu,\mu'$ such that
\[
t^*_{\mu,\mu',n} =  t_{[i],n} \asymp \frac{1}{2}\frac{\log(n/n_0)}{\log(1/\lambda_{[i]})}\, ,
\]
for some initial sample complexity $n_0$, i.e. the sample complexity to distinguish $\mu$ and $\mu^\prime$ without the action of $P$.
\end{theorem}
\begin{figure}[h!]
\label{FIG:markov_phase}
    \begin{center}
      \includegraphics[width= 0.8\textwidth]{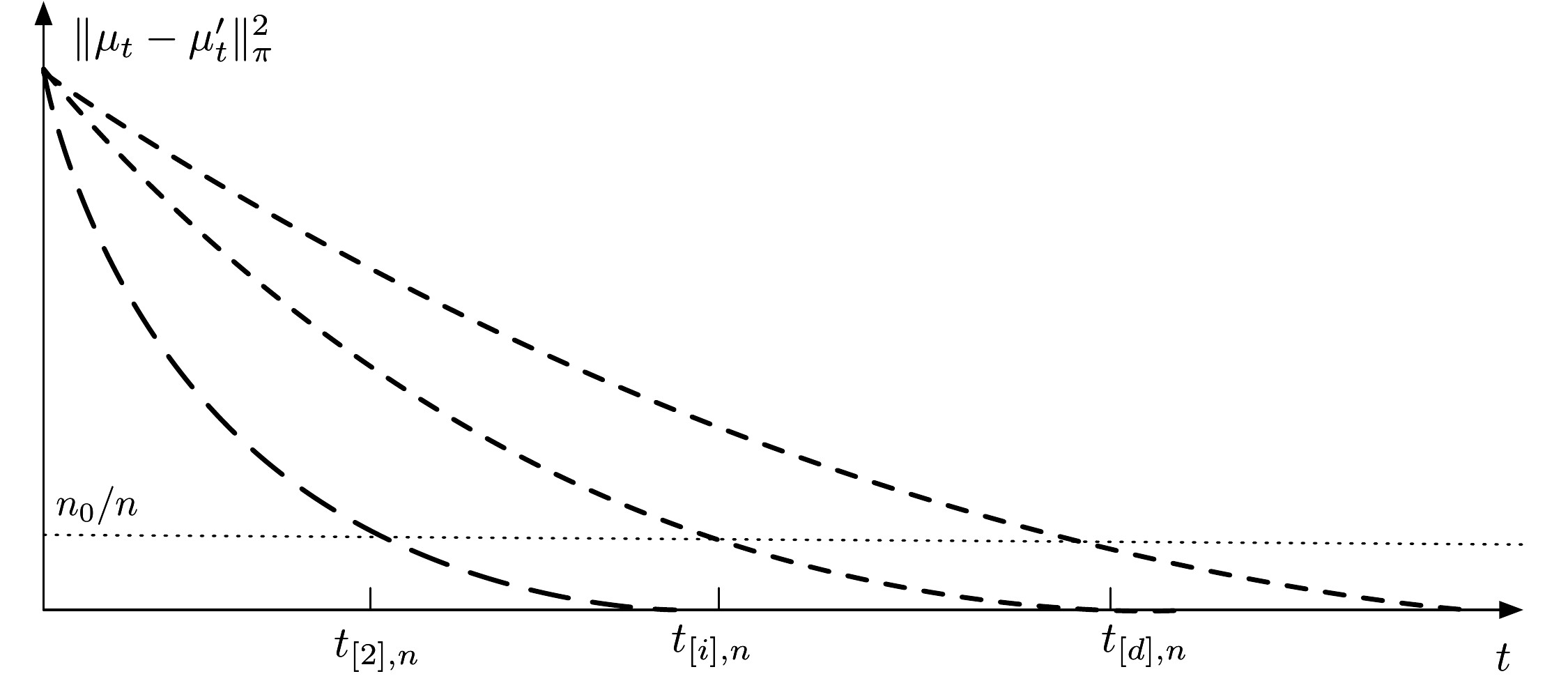}
    		\caption{For different cases of starting distributions, the time at which the squared distance between the two distributions is much smaller than $n_0/n$, and the hypothesis testing problem becomes statistically impossible, can vary greatly.}
    \end{center}
\end{figure}

Yet again, this establishes that there is not only one important time, such as mixing time, describing the loss of information in this problem, but that this can happen at different timescales; alternatively, the loss of information is not only driven by the spectral gap, but may depend on all the eigenvalues. 

\section{Applications}
\label{sec:applications}

In this section, we consider several concrete examples of Markov Chains on
commonly studied graph topologies. $P$ denotes the transition matrix of the
Markov Chain in the corresponding section. The ratio of the second largest
eigenvalue to the smallest, when ordered by absolute value, gives the
\emph{statistical window} in each case. The corresponding eigenvectors shed
light on the pairs of starting distributions that give rise to the extreme ends
of this window. 

Let $u_{[2]}$ and $u_{[d]}$ denote the left eigenvectors corresponding to the
second largest and smallest eigenvalues when ordered by absolute values, i.e.
$\lambda_{[2]}$ and $\lambda_{[d]}$, respectively. For a sufficiently small
$\alpha$, we can consider the pairs of distributions given by $(\mu,
\mu^\prime)$ and $(\gamma, \gamma^\prime)$, where $\mu = \pi + \alpha u_{[2]}$,
$\mu^\prime = \pi - \alpha u_{[2]}$, $\gamma = \pi + \alpha u_{[d]}$ and
$\gamma^\prime = \pi - \alpha u_{[d]}$; $\alpha$ is chosen to be small enough so that $\mu, \mu^\prime, \gamma, \gamma^\prime$ are all valid probability distributions and 
$\mu, \mu^\prime, \pi$ pairwise satisfy the $\varepsilon$-bounded likelihood
condition, as do the distributions $\gamma, \gamma^\prime, \pi$. Note that by
Theorems~\ref{THM:main1} and~\ref{THM:main2}, the sample complexity of
distinguishing between the pairs $(\mu, \mu^\prime)$ and $(\gamma,
\gamma^\prime)$ is the same up to factors depending only on $\varepsilon$ and
$\delta$---in each case it is given by $1/4\alpha^2 = \lnormpi{\mu -
\mu^\prime}^{-2} = \lnormpi{\gamma - \gamma^\prime}^{-2}$. However, when
considering the statistical complexity of distinguishing between $(\mu_t,
\mu^\prime_t)$, the sample complexity grows as $\lambda_{[2]}^{-2t}$, whereas
for distinguishing between $(\gamma_t, \gamma^\prime_t)$ it grows as
$\lambda_{[d]}^{-2t}$. Thus, despite having roughly the same sample complexity
at time $t = 0$, at a later time the ratio of sample complexities grows to be
as large as $(\lambda_{[2]}/\lambda_{[d]})^{2t}$. In each of the concrete
examples below, we discuss lower bounds on $\lambda_{[2]}$ and upper bounds on
$\lambda_{[d]}$ and the ratio of these bounds gives a lower bound on the
statistical window. In some of the examples, we explicitly derive the
eigenvectors associated with these eigenvalues and discuss the associated pairs
of initial distributions.

\subsection{Random Walk on a Bipartite Clique}

We start with the simple example of a random walk on a bipartite clique. Let
$\cX = L \cup R$ be an equi-partition and let $E = \{\{x, x^\prime\} | x \in L,
x^\prime \in R \}$, then the transition matrix $P$ is given by: 
\begin{align*}
	\fP[X_s = x| X_{s-1} = x^\prime] &= {2}/{d} \,\quad\quad\text{if } \{x , x^\prime\} \in E
\end{align*}

The transition matrix $P$ has exactly two non-zero eigenvalues, $1$ and $-1$;
thus $\lambda_{[2]} = -1$ and $\lambda_{[d]} = 0$. The eigenvector
corresponding to the eigenvalue $-1$ has negative entries on one side of the
bi-partite graph and positive entries on the other. Thus, the component of this
eigenvector in the distribution controls the imbalance of the distribution between the two sides.  If
the initial distributions satisfy $\mu(L) = \mu^\prime(L)$, they become
indistinguishable after just one step of the Markov Chain. On the other hand,
the difference $|\mu(L) - \mu^\prime(L)|$ remains unaffected by the Markov
chain. Thus, the problem at any time $t \geq 1$ remains exactly as hard as the
problem at time $t = 1$, i.e. all initial information except for the
\emph{starting side} is lost in exactly one time step. So the \emph{statistical
window} is infinite in this case, for $t \geq 1$.

\subsection{Random Walk on the Cycle} 
\begin{definition}[Random Walk on the $d$-Cycle]
	Let $\cX = \{0, 1, \ldots, d - 1\}$ be the $d$ nodes of the cycle. Let $P$
	be the Markov chain on $\cX$, where,
	\begin{align*}
		\fP(X_s = i | X_{s-1}=j) &=
		\begin{cases}
			\frac{1}{2} & \text{if $i \equiv j \pm 1 (\modulo d)$} \\
			0 & \text{otherwise}
		\end{cases}	
	\end{align*}
\end{definition}
The spectral properties of $P$ are well known (see e.g.~\cite[Chap.
12.3]{LPW:2008}); we summarize them in the following lemma.
\begin{lemma}
	For any $d \geq 3$, the eigenvalues of $P$ are given by
	$\cos({2 \pi i}/{d})$ for $i \in \{0, \ldots, d - 1\}$; the
	(right and left) eigenvector $u_i = (u_{i,0}, \ldots, u_{i, d-1})$
	corresponding to eigenvalue $\cos({2 \pi i}/{d})$, is given
	by $u_{i,k} = \cos({2 \pi i k}/{d})$.
\end{lemma}
Let us first consider a cycle of length $d$ with $d \equiv 0 (\modulo 4)$. In
this case, $\lambda_{[2]} = -1$ and $\lambda_{[d]} = 0$. The associated
eigenvectors, $u_{[2]}$ and $u_{[d]}$, are (up to scaling) given by:
$u_{[2]}(i) = 1$ for even $i$ and $-1$ for odd $i$; $u_{[d]}(i) = 1$ for $i
\equiv 0 (\modulo 4)$, $-1$ for $i \equiv 2 (\modulo 4)$ and $0$ for $i \equiv
1 (\modulo 2)$. Now consider the pairs $\mu, \mu^\prime$, where $\mu = \pi +
\alpha v_{[2]}$ and $\mu^\prime = \pi - \alpha v_{[2]}$, and $\gamma,
\gamma^\prime$, where $\gamma = \pi + \alpha v_{[d]}$ and $\gamma^\prime = \pi
- \alpha v_{[d]}$. In the first case, it is easy to see that the probability
mass on \emph{odd} and \emph{even} nodes is noticeably different under $\mu$
and $\mu^\prime$, and this will remain so in perpetuity. On the other hand,
starting from $\gamma$ or $\gamma^\prime$, stationarity is achieved in one
step. Observe that in this case $n^*_{\mu, \mu^\prime, 0} \asymp n^*_{\gamma,
\gamma^\prime, 0}$, so initially the two problems are roughly equally hard;
however, as $t$ increases (in the simple case of cycle lengths being multiples of 4, even for $t = 1$) the difference between the statistical hardness of
these problems differs dramatically. This behavior will be approximately
replicated with cycles of any length provided $d$ is large enough; in
particular we always have $|\lambda_{[2]}| = 1 - O(1/d^2)$ and $\lambda_{[d]} =
O(1/d)$, and so the statistical window is of size $d^{2t}$.

\subsection{Random Walk on the Line Graph}

The random walk on the line graph with $d$ nodes is very similar to that on the cycle. In fact, the walk can be viewed as a projection of the random walk on a cycle with $2(d - 1)$ nodes. 

\begin{definition}[Random Walk on $d$-Line] Let $\cX = \{0, \ldots, d - 1\}$ be the $d$ nodes on a line. Let $P$ be the Markov Chain on $\cX$, where,
	\[ \fP(X_s = i| X_{s-1} = j) = 
	\begin{cases}
		\frac{1}{2} & \text{for } i \in {j - 1, j + 1} , j \in \{1, \ldots, d - 2 \} \\
		1 & \text{if } j = 0 \wedge i = 1, \text{ or } j = d - 1 \wedge i = d - 2 \\
		0 & \text{otherwise}
	\end{cases}
	\]
\end{definition}

As in the case of the cycle, the spectrum is explicitly known (cf.~\cite[Chap.
12.3]{LPW:2008}). The result is stated as the following lemma; this implies a statistical window of at least $d^{2t}$.

\begin{lemma}
	For any $d \geq 3$, the eigenvalues of $P$ are given by
	$\cos({\pi i}/{(d - 1)})$ for $i \in \{0, \ldots, d - 1\}$; the
	right eigenvector $u_i = (u_{i,0}, \ldots, u_{i, d-1})$
	corresponding to eigenvalue $\cos({\pi i}/{(d - 1)})$, is given
	by $u_{i,k} = \cos({\pi i k}/{(d - 1)})$.
\end{lemma}

\subsection{Random Walk on the Regular Block Model}

We use a variant of the stochastic blockmodel~\citep{HolLasLei83} where the graph is regular (as
opposed to approximately regular). Note that the model is completely deterministic.

\begin{definition}[Regular Blockmodel]
	A regular blockmodel with $k$ blocks on $d$ nodes with degrees $(\Delta_{i,j})$ for $1
	\leq i, j \leq k$ and $\Delta_{i,j} = \Delta_{j, i}$ is defined as follows:
	The vertex set $V$ is partitioned as $V =  V_1  \cup V_2 \cup \cdots \cup
	V_k$, with $|V_i| = d/k$. The induced subgraph $G_i = (V_i, E(V_i))$ is a
	$\Delta_i$ regular graph for each $i$, and the subgraph $G_{i, j} = (V_i
	\cup V_j, E(V_i \cup V_j) \setminus (E(V_i) \cup E(V_j)))$ is a
	$\Delta_{i,j}$ regular bipartite graph for all $i, j$, $i \neq j$.
\end{definition}

\begin{proposition}
	\label{prop:regular-blockmodel}
	There exist regular block models with $k = 2$ blocks on $d$ nodes, satisfying
	\[ 
		|\lambda_{[2]}(P)| = 1 - o(1)\, , \quad \text{and} \quad |\lambda_{[d]}(P)| \asymp {1}/{d} 
	\]
\end{proposition}

For blockmodels with $k = 2$, the eigenvector corresponding to $\lambda_{[2]} =
\lambda_2$, correlates strongly (in fact for the regular blockmodel with equal
sized blocks, exactly), with the block structure. Thus, if $\mu$ and $\mu^\prime$ start
off with significantly different probability mass on the two blocks, the
statistical problem remains easy essentially until mixing time. On the other
hand, if they have the same probability mass on the individual blocks (even
though the distributions may differ on the blocks significantly), in typical cases,
the statistical problem becomes hard quickly, e.g. if each block is an
expander.

\begin{remark} 
	A special case of regular graphs is the class of Ramanujan graphs~\citep{LubPhiSar88}; the
	eigenvalues of the transition matrix of a random walk on a Ramanujan graph
	are $\pm 1$ or satisfy $|\lambda| = O(1/\sqrt{d})$. In the non-bipartite
	case, there is no guarantee that the statistical window is large.
\end{remark}

\subsection{Pachinko random walk}
\label{sec:tree}

We introduce the following random walk inspired by the Japanese pinball game of {\em Pachinko}, on $\cX = [d]$ where $d=2^r$ and the space $\cX$ is understood as the leaves of a dyadic tree of height $r$. It allows to further illustrate the statistical window phenomenon.

\begin{definition}[Pachinko Random Walk]
Let $P$ be the Markov chain on the $d=2^r$ leaves of a dyadic tree such that for two leaves $i$ and $j$ with first common ancestor at height $\ell$ between $0$ and $r$, we have
\[
\fP(X_{s}=i | X_{s-1}=j) = p_\ell = {\beta_\ell}/{2^{\ell-1}}\, ,
\]
where $\beta_0>\ldots>\beta_r$ are positive real numbers that sum to 1.
\end{definition}

\begin{figure}[h!]
\label{FIG:pachinko}
    \begin{center}
      \includegraphics[width= 0.8\textwidth]{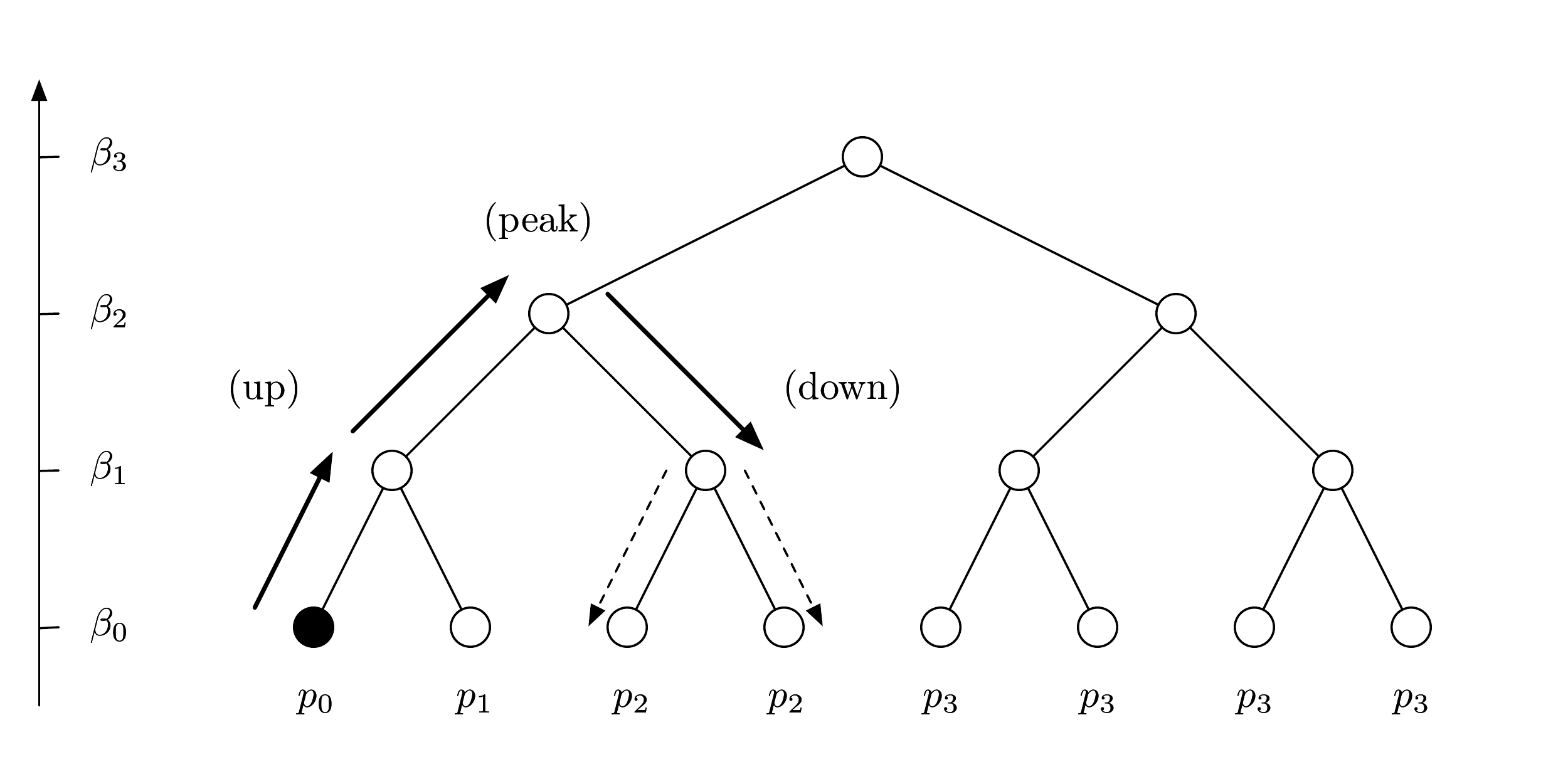}
    		\caption{Analogously to the game of Pachinko, a ball starts at a leaf, goes 	{\em (up)} in the dyadic tree, {\em (peaks)} with probability $\beta_\ell$ at height $\ell$, and goes {\em (down)} on the other side of the highest point, going left or right 	independently with probability $1/2$ at each further descendant node on the way down, thus stopping uniformly at random in one of the $2^{\ell-1}$ leaves. In this figure $r=3$, and we represent a trajectory that peaks at height $\ell=2$.}
    \end{center}
\end{figure}

This can be understood as a random walk on a graph, with a large amount of structure, that can be a consequence of an underlying geometry: at every level, each half of the vertices is ``very far'' from the other half, and jumping from one half to the other is less probable than staying in the same half.

\begin{proposition}
	\label{prop:tree}
For each $k$ between $2$ and $r+1$, there exists an eigenvalue $\gamma_k$ of multiplicity $2^{k-2}$, associated to the height $\ell = r+2-k$. It is given by $\gamma_k = \beta_0+\ldots+\beta_{r+1-k} - \beta_{r+2-k} = \lambda_i$, for $2^{k-2}+1\le i \le 2^{k-1}$. The $2^{k-2}$ vectors associated to the $2^{k-2}$ nodes at height $\ell = r+2-k$ - with coefficients equal to $1$ for their left descendants, $-1$ for their right descendants, and 0 otherwise - are eigenvectors with eigenvalue $\gamma_k$. For $k=1$, $\gamma_1 = \lambda_1=1$ is an eigenvalue with associated simplex eigenvector $\pi = \mathbf{1}/d$. The eigenvalues satisfy $\gamma_1 >\gamma_2>\ldots>\gamma_{r+1}>0$. 
\end{proposition}
This situation is summarized in the case $r=3$ of Figure~\ref{FIG:pachinko} in the following table.
\begin{center}
    \begin{tabular}{ | l | l | l | p{5.5cm} |}
    \hline
        $k$ & $\ell$ & $\gamma_k$ & $\|.\|_\pi$ orthonormal eigenvectors \\ \hline
    $k=1$ &  & $\gamma_1= \beta_0+\beta_1+\beta_2+\beta_3=\lambda_1 =1$ & \makecell[l]{$v_1 = (1,1,1,1, 1,1,1,1)/8$} \\ \hline
    $k=2$ & $\ell=3$ & $\gamma_2= \beta_0+\beta_1+\beta_2-\beta_3=\lambda_2 $ & \makecell[l]{$v_2 = (1,1,1,1, -1,-1,-1,-1)/8$} \\ \hline
    $k=3$ & $\ell=2$ & $\gamma_3= \beta_0+\beta_1-\beta_2=\lambda_3 = \lambda_4 $ & \makecell[l]{$v_3 = (1,1,-1,-1,0,0,0,0)/2^{5/2} $ \\
    $v_4 = (0,0,0,0,1,1,-1,-1)/2^{5/2} $}\\ \hline
    $k=4$ & $\ell=1$ & $\gamma_4= \beta_0-\beta_1=\lambda_5=\lambda_6=\lambda_7=\lambda_8 $ & \makecell[l]{$v_5 = (1,-1,0,0,0,0,0,0)/4 $\\$v_6 = (0,0,1,-1,0,0,0,0)/4$, etc.}\\
    \hline
    \end{tabular}
\end{center}

As a consequence, if between height $\ell=r+1-k$ and $\ell+1=r+2-k$, there is a large gap between the probabilities $\beta_\ell$ and $\beta_{\ell+1}$ (i.e., it is much harder for a particle to ``jump'' to height $\ell+1$ than to height $\ell$), we have that
\[
\gamma_{k} - \gamma_{k+1} = 2\beta_{\ell}-\beta_{\ell+1}>\beta_\ell-\beta_{\ell+1}>0\, ,
\]
which implies a large gap between the eigenvalues of eigenvectors associated to height superior to $\ell$ and those associated to a height less or equal to $\ell$. From a statistical point of view, in light of Theorem~\ref{THM:main1} and similarly to Theorem~\ref{THM:window}, this implies that difference between $\mu$ and $\mu'$ which is observable at height $\ell+1$ or above, i.e. difference of mass between two sides of nodes at these height) will be statistically observable for much larger $t$ than difference at lower levels. 

\subsection{Random Markov chains}
The notion of random Markov chains, and in general of random walks in random media, has been thoroughly studied. We consider here the case of random reversible Markov chains, as studied by \cite{BCC:2010spectrum-two-examples, BCC:2011spectrum}. 

\begin{definition}
	\label{defn:random-markov}
	Consider a finite connected undirected graph $G = (\cX, E)$; for every $\{i,j\} \in E$,
	let $U_{i,j}$ be drawn in an i.i.d. manner from a distribution on the
	positive part of the real line with bounded second moments; we set $U_{i, j}
	= U_{j, i}$. (An edge $\{i \}$ is allowed, which corresponds to a self-loop
	at $i$.) The coefficients of $P$ are obtained by normalization to a
	stochastic matrix
	\[
		P_{i,j} = U_{i,j}/\sum_{x \in \cX} U_{i,x}\, .
	\]
\end{definition}

Understanding the spectra of random (symmetric) matrices has been intensely
studied in recent years (see
e.g.~\cite{TV:2010random,Cha:2009singular,RudVer08,Ver:2014invertibility}). Below, we
use results from this literature to understand the behavior of
$\lambda_{[2]}(P)$ and $\lambda_{[d]}(P)$; again, this yields a statistical window of $d^{2t}$.

\begin{proposition}
	\label{prop:random-markov}
	Let $P$ be the transition matrix of a random Markov chain as defined in
	Defn.~\ref{defn:random-markov} with $G(\cX, E)$ being the complete graph. As $d \rightarrow \infty$
	$d$, it holds with probability going to $1$ that
	\[
		|\lambda_{[2]}(P)| \asymp {1}/{\sqrt{d}}\, , \quad \text{and} \quad |\lambda_{[d]}(P)| \asymp {1}/{d^{3/2}}\, ,
	\]
\end{proposition}

\subsection{Random Walk on the Hypercube}

The hypercube on $d = 2^k$ nodes can be viewed as a graph with the node set denoted by $\{-1, 1\}^{k}$. We first consider the standard random walk on the hypercube which is defined below; for $x, x^\prime \in \{-1, 1\}^{k}$, let $|x - x^\prime|_H := \frac{1}{2} \Vert x - x^\prime \Vert_1$ denote the Hamming distance between $x$ and $x^\prime$. 

\begin{definition}
	Let $\cX = \{-1, 1\}^{k}$ be the $d$ nodes of the hypercube. Let $P$ be the Markov chain on $\cX$, where
	\begin{align*}
		\fP(X_s = x | X_{s - 1} = x^\prime) &= 
		\begin{cases}
			\frac{1}{d}  & \text{if } |x - x^\prime|_H = 1 \\
			0 & \text{otherwise}
		\end{cases}
	\end{align*}
\end{definition}

The spectral properties of $P$ are summarized by the following lemma
\citep[c.f.][Chap 12.4]{LPW:2008}; this implies a statistical window of at
least $d^{2t}$.

\begin{lemma}
	For $d = 2^k$ (with $k \geq 1$), the eigenvalues of $P$ are given by $1 -
	{2j}/{k}$, $j = 0, \ldots, k$. The eigenvalue $1 - {2j}/{k}$ appears
	with multiplicity ${k \choose j }$; the corresponding eigenvectors $u$ are
	given by considering sets of size $j$. A set $S \subseteq \{1, \ldots, k
	\}$, $|S| = j$, yields the eigenvector $u(x) = \prod_{i \in S} x_i$, for $x
	\in \{-1, 1\}^k$. (These are the so-called parity functions.)
\end{lemma}

\paragraph{Product Distributions on the Hypercube.} 

One can consider other product distributions on the hypercube. Consider a two
state Markov chain, with states denoted by $\cX_2 = \{ -1, 1\}$ and transitions
given as:

\begin{definition}
	\label{defn:two-state-chain}
	Let $\cX_2 = \{-1, 1 \}$, $0 < p, q \leq 1$, define $P : \cX_2 \times \cX_2
	\rightarrow [0, 1]$ as $P(-1, -1) = 1 - p$, $P(-1, 1) = p$, $P(1, -1) = q$,
	$P(1, 1) = 1 - q$. The Markov chain on $\cX_2$ defined by $P$ is given by,
	\begin{align*}
		\fP(X_s = x |  X_{s - 1 }= x^\prime ) = P(x, x^\prime) 
	\end{align*}
\end{definition}

The eigenvectors of $P$ are easiest to express as functions from $\cX_2
\rightarrow \reals$. The right eigenvector corresponding to eigenvalue $1$ is
given by $u(x) = 1$ for $x \in \cX_2$. The stationary distribution is given by
$\pi(-1) = q/(p + q)$ and $\pi(1)= p/(p + q)$, let $\xi = \E_{x \sim \pi} [x] =
(p - q)/(p + q)$. The second eigenvalue is $1 - (p + q)$ and the corresponding
eigenvector is given by $u(x) = (x - \xi)/\sqrt{1 - \xi^2}$. 

\begin{definition}[Product Chain on the Hypercube] 
	Let $d = 2^k$, $\cX = \{-1, 1\}^k = \cX_2 \otimes \cdots \otimes \cX_2$, let
	$P^{(1)}, \ldots, P^{(k)}$ be transition matrices of the chain on $\cX_2$
	defined in Defn.~\ref{defn:two-state-chain} with parameters $(p^{(1)},
	q^{(1)}), \ldots, (p^{(k)}, q^{(k)})$, and let $w_1, \ldots, w_k$ be
	positive weights such that $\sum_{i} w_k = 1$. Then, for $x, x^\prime \in
	\cX$, we have the following Markov chain:
	\begin{align*}
		\fP(X_s = x | X_{s-1} = x^\prime) = 
		\begin{cases} 
			\sum_{j = 1}^k w_j P^{(j)}(x_j, x^\prime_j) & \text{if $| x - x^\prime|_H = 1$} \\
			0 & \text{otherwise}
		\end{cases}
	\end{align*}
\end{definition}

The eigenvectors and eigenvalues of the product Markov chain  on the hypercube
are easily defined through the eigenvectors and eigenvalues of Markov chain
defined on $\cX_2$. The following lemma follows from results stated
in~\cite[Chap 12.4]{LPW:2008}.
\begin{lemma}
	Let $P$ be the transition matrix of the product Markov chain obtained using the transition matrices $P^{(1)}, \ldots, P^{(k)}$ of chains on $\cX_2$. Let $u^{(i)}$ denote the eigenvector of $P^{(i)}$ with eigenvalue $1 - (p^{(i)} + q^{(i)})$, then for each subset $S \subseteq \{1, \ldots, k\}$, define $u_S : \cX \rightarrow \reals$ as follows:
	\begin{align*}
		u_S(x) = \prod_{i \in S} u^{(i)}(x_i) = \prod_{i \in S} \frac{x_i  - \xi^{(i)}}{\sqrt{ 1 - (\xi^{(i)})^2}}
	\end{align*}
	where $\xi^{(i)} = (p^{(i)} - q^{(i)})/(p^{(i)} + q^{(i)})$. Then $u_S$ is an eigenvector of $P$ with eigenvalue $1 - \sum_{i \in S} w_i (p^{(i)} + q^{(i)})$.
\end{lemma}

\begin{remark} It is easily observed that if we set $p^{(i)} = q^{(i)} = 1$ in all the chains and $w_i = \frac{1}{k}$ for each $i$, then we get exactly the standard random walk on the hypercube with $d = 2^k$ vertices.
\end{remark}

\bibliographystyle{plainnat}
\bibliography{statebib2}

\appendix

\label{app:info_theory}
\section{Distances and divergences between distributions}
\label{APP:notation}

We recall the following notions, for two distributions $\mu,\mu'$ on a finite set $\cX$. The total variation distance, denoted by $d_{\sf TV}$, is defined as
\[
d_{\sf TV}(\mu,\mu') = \frac{1}{2}\sum_{x \in \cX}|\mu_x - \mu'_x| = \max_{X \subseteq \cX} |\mu(X) - \mu'(X)| \, .
\]
The Kullback--Leibler divergence, denoted by $\KL$, is defined as
\[
\KL(\mu,\mu') = \sum_{x \in \cX}\mu_x \log\left(\frac{\mu_x}{\mu'_x} \right) \, .
\]
The chi-square divergence, denoted by $D_{\chi^2}$, is defined as
\[
D_{\chi^2}(\mu,\mu') = \sum_{x \in \cX}\mu'_x \left(\frac{\mu_x}{\mu'_x} -1\right)^2 \, .
\]
The Hellinger distance, denoted by $H$, is defined by
\[
H^2(\mu,\mu') = \sum_{x \in \cX}\mu'_x \left(\sqrt{\frac{\mu_x}{\mu'_x}} -1\right)^2 \, .
\]
We will also make frequent use of the following standard inequalities (see e.g.~\cite{Tsy08}). 
\begin{align*}
	\frac{1}{2} H^2(\mu, \mu^\prime) &\leq d_{\sf TV}(\mu, \mu^\prime) \leq H(\mu, \mu^\prime) \\
	d_{\sf TV}(\mu, \mu^\prime) &\leq \sqrt{\frac{1}{2} \KL(\mu, \mu^\prime)} \tag{Pinsker's Inequality} \\
\end{align*}

%
%
\section{Proofs}
\label{sec:app}
\subsection{Proofs from Section~\ref{sec:samplecomp}}

\begin{proof}[of Theorem~\ref{THM:main1}]
Recalling that $\mu_t = \mu P^t$, the following classical result holds, for  $\mu$ and $\mu^\prime$ be distributions over $\cX$. 

\[	%
\mu = \pi + \sum_{i = 2}^d \ippi{u_i}{\mu} u_i \quad \text{and   } \mu^\prime = \pi + \sum_{i = 2}^d \ippi{u_i}{\mu^\prime} u_i\, ,
\]

Indeed, $\mu$ and $\mu^\prime$ are treated as vectors in $\reals^d$ and can be
	expressed in terms of the basis vectors $u_1, \ldots, u_d$. The coefficients
	are given by $\alpha_i = \ippi{\mu}{u_i}$ (resp. $\alpha_i^\prime  =
	\ippi{\mu^\prime}{u_i}$). As $\mu$ and $\mu^\prime$ are distributions,
	Lemma~\ref{lem:markov-basis} give us that $\alpha_1 = \alpha_1^\prime = 1$. 

	As the $u_i$'s are eigenvectors, we have $\mu_t = \mu P^t = \pi + \sum_{i =
	2}^d \alpha_i \lambda^t_i u_i$ and $\mu^\prime_t = \mu^\prime P^t = \pi +
	\sum_{i = 2}^d \alpha^\prime_i \lambda^t_i u_i$. Then, by using the orthonormality of $u_i$ with respect to $\ippi{\cdot}{\cdot}$, we have:

	\[
\lnormpi{\mu_t - \mu_t^\prime}^2 = \sum_{i = 2}^d \lambda_i^{2t} (\ippi{u_i}{\mu} - \ippi{u_i}{\mu^\prime})^2
\]

This norm between $\mu_t$ and $\mu'_t$ can be compared to other notions of distances between distributions. In particular, this can be done for the Hellinger distance between $\mu_t$ and $\mu'_t$. Let $\mu_{t,x}/\mu^\prime_{t,x} =: 1 + \gamma_{t,x}$. 
	 It then holds that
	\begin{align*}
		H^2(\mu_t,\mu'_t) &=  \sum_x \mu^\prime_{t,x} (1 -\sqrt{1+\gamma_{t,
		x}})^2 \\
		&= 2 \sum_x \mu'_{t,x} \Big(1+\frac{\gamma_{t,x}}{2} -
		\sqrt{1+\gamma_{t,x}} \Big)\\
		&\ge \frac{\varepsilon^{3/2}}{8} \sum_x \mu'_{t,x} \gamma_{t,x}^2\, \quad
		\text{(see below)}\\
		&= \frac{\varepsilon^{3/2}}{8} \sum_x \frac{(\mu_{t,x} -
		\mu_{t,x}^\prime)^2}{\mu^\prime_{t,x}} 
\end{align*}

The first inequality can be justified in the following way: Define the function
$f_c$ on $[-1,\infty)$ by 
\[
	f_c(t)=1+t/2-ct^2-\sqrt{1+t}\, .
\] 
Basic computations yield that for all positive $c$, $f_c(0)=0$ and $f'_c(0)=0$.
Furthermore, we have that 
\[
	f''_c(t) = -2c+\frac{1}{4(1+t)^{3/2}}\, .
\]
As a consequence, $f_c$ is convex in $t\in [\varepsilon-1,1/\varepsilon-1]$ for $c=\varepsilon^{3/2}/8$, so $f_c(t) \ge 0$ on this interval. We therefore have that $1+\frac{\gamma_{t,x}}{2} - \sqrt{1+\gamma_{t,x}} \ge \frac{\varepsilon^{3/2}}{8} \gamma_{t,x}^2$ for $\mu_{t,x}/\mu'_{t,x} \in [\varepsilon,1/\varepsilon]$. As a consequence, it holds that
\begin{align*}
	H^2(\mu_t,\mu'_t) &\ge \frac{\varepsilon^{5/2}}{8} \sum_x \frac{(\mu_{t,x} - \mu_{t,x}^\prime)^2}{\pi_x}\\
&= \frac{\varepsilon^{5/2}}{8} \|\mu_t - \mu'_t\|_\pi^2
\end{align*}
We use this property to control directly the probability of error of the likelihood-ratio test $\psi_{\sf LR}$
\begin{align*}
\max_{\nu \in \{\mu,\mu'\}} \fP^{\otimes n}_{\nu_t}(\psi_{\sf LR} \neq \nu)&\le \fP^{\otimes n}_{\mu_t}(\psi_{\sf LR} \neq \mu) + \fP^{\otimes n}_{\mu'_t}(\psi_{\sf LR} \neq \mu')\\
&= 1-  \big(\fP^{\otimes n}_{\mu'_t}(\psi_{\sf LR} \neq \mu) - \fP^{\otimes n}_{\mu_t}(\psi_{\sf LR} \neq \mu)\big)\\
&=1-d_{\sf TV}(\mu_t^{\otimes n},\mu_t^{\prime \otimes n})
\end{align*}
Indeed, one of the event $X$ realizing the total variation distance between the two distributions is the one on which $\mu_t^{\prime \otimes n}\big((X_i)_{i \in [n]} \big)$ is greater or equal than $\mu_t^{\otimes n}\big((X_i)_{i \in [n]} \big)$, i.e. where the output of $\psi_{\sf LR}$ is $\mu'$.
We have, by properties of the Hellinger distance, that
\[
1-d_{\sf TV}(\mu^{\otimes n}_t,\mu^{\prime \otimes n}_t) \le 1-\frac{1}{2}H^2(\mu_t^{\otimes n},\mu_t^{\prime \otimes n}) = \Big(1-\frac{1}{2} H^2(\mu_t,\mu_t') \Big)^n\, .
\]
The last equality allows, by tensorization, to relate directly this probability of error to a quantity depending separately on $(\mu_t,\mu'_t)$, and $n$. As a consequence, we have
\[
\max_{\nu \in \{\mu,\mu'\}} \fP^{\otimes n}_{\nu_t}(\psi_{\sf LR} \neq \nu) \le \Big(1-\frac{\varepsilon^{5/2}}{16} \|\mu_t - \mu'_t\|_\pi^2 \Big)^n \le e^{-n\frac{\varepsilon^{5/2}}{16} \|\mu_t - \mu'_t\|_\pi^2}
\]
As a consequence, for $n\ge 16\varepsilon^{-5/2} \log(1/\delta)/ \|\mu_t - \mu'_t\|_\pi^2$, the probability of error is indeed less than $\delta$.
\end{proof}

\begin{proof}[of Theorem~\ref{THM:main2}]
We recall from the proof of Theorem~\ref{THM:main1} that
\[
\lnormpi{\mu_t - \mu_t^\prime}^2 = \sum_{i = 2}^d \lambda_i^{2t} (\ippi{u_i}{\mu} - \ippi{u_i}{\mu^\prime})^2\, .
\]
We compare this distance to the Kullback--Leibler divergence $\KL(\mu_t,\mu'_t)$. Let $\mu_{t,x}/\mu^\prime_{t,x} =: 1 + \gamma_{t,x}$; as $\sum_x \mu_{t,x} = \sum_x
	\mu^\prime_{t,x}$, it follows that $\sum_x \mu^\prime_{t,x} \gamma_{t,x} = 0$.  Then
	consider the following:

	\begin{align*}
		\KL(\mu_t, \mu_t^\prime) &= \sum_{x} \mu_{t,x} \ln\frac{\mu_{t,x}}{\mu^\prime_{t,x}} = \sum_x \mu^\prime_{t,x} (1 + \gamma_{t,x}) \ln(1 + \gamma_{t,x}) \\
		&\leq \sum_x \mu^\prime_{t,x} (1 + \gamma_{t,x}) \gamma_{t,x} &\text{Using the fact that $\ln(1 + t) \leq t$} \\
		&=\sum_x \mu^\prime_{t,x} \gamma_{t,x}^2 &\text{As $\sum_x \mu_x^\prime \gamma_x = 0$} \\
		&= \sum_x \frac{(\mu_{t,x} - \mu_{t,x}^\prime)^2}{\mu^\prime_{t,x}} \\
		&\le \frac{1}{\varepsilon} \sum_x \frac{(\mu_{t,x} - \mu_{t,x}^\prime)^2}{\pi_{x}} &\text{As $\frac{1}{\mu^\prime_{t, x}} \leq \frac{1}{\varepsilon \pi_x}$} \\
		&= \frac{1}{\varepsilon}\|\mu_t-\mu'_t\|_\pi^2\,.
	\end{align*}

To give a lower bound on the probability of error, we have
\[
\inf_{\psi} \max_{\nu \in \mu, \mu'} \fP^{\otimes n}_{\nu_t}(\psi \neq \nu) \ge
\inf_{\psi}
	\frac{1}{2}(\fP^{\otimes n}_{\mu_t}(\psi \neq \mu)+\fP^{\otimes
	n}_{\mu'_t}(\psi \neq \mu')) \ge \frac{1-d_{\sf TV}(\mu_t^{\otimes
	n},\mu_t^{\prime \otimes n})}{2}
\]
The above holds by using the definition of total variation distance as the supremum of the difference in probability for all events, and using the event $\psi=\mu$, with an infimum taken over all tests $\psi$. Furthermore, by Pinsker's inequality and by the tensorization properties of the Kullback--Leibler divergence, we have that
\[
d_{\sf TV}(\mu_t^{\otimes n},\mu_t^{\prime \otimes n}) \le  \sqrt{\KL(\mu^{\otimes n}_t,\mu^{\prime \otimes n}_t)/2} = \sqrt{n \KL(\mu_t,\mu^{\prime}_t)/2}\, .
\]
As a consequence, it holds that
\[
\inf_{\psi} \max_{\nu \in \mu, \mu'} \fP^{\otimes n}_{\nu_t}(\psi \neq \nu) \ge \frac{1}{2} - \frac{1}{2}\sqrt{\frac{n}{2 \varepsilon} \|\mu_t - \mu'_t\|_\pi^2}\, .
\]
For any $n\le 8\varepsilon \delta^2/\|\mu_t-\mu'_t\|_\pi^2$, the probability of error is at least $1/2-\delta$.
\end{proof}

\subsection{Proofs from Section~\ref{SEC:general}}

\begin{proof}[of Theorem~\ref{THM:general}]
As in the proof of Theorem~\ref{THM:main1}, it holds that

\[
\max_{\nu \in \{\mu,\mu'\}} \fP^{\otimes n}_{\nu_t}(\psi_{\sf LR} \neq \nu) \le 1-d_{\sf TV}(\mu^{\otimes n}_t,\mu^{\prime \otimes n}_t) \le 1-\frac{1}{2}H^2(\mu_t^{\otimes n},\mu_t^{\prime \otimes n}) = \Big(1-\frac{1}{2} H^2(\mu_t,\mu_t') \Big)^n\, .
\]
Furthermore, we have that
\begin{align*}
H^2(\tilde \mu_t,\tilde \mu'_t) &=\sum_{x}\big(\sqrt{\tilde \mu_{t,x}}-\sqrt{\tilde \mu'_{t,x}}\big)^2\\
&=\sum_{x}(\sqrt{(1-\eta)\mu_{t,x} + \eta \beta_{t,x}}-\sqrt{(1-\eta)\tilde \mu'_{t,x}+\eta \beta_{t,x}})^2\\
&\le (1-\eta) \sum_{x}\big(\sqrt{\mu_{t,x}}-\sqrt{\mu'_{t,x}}\big)^2 \\
&= (1-\eta) H^2(\mu_t,\mu'_t)\, .
\end{align*}
Indeed the inequality is a consequence of $\sqrt{a+c}-\sqrt{b+c} \le \sqrt{a}-\sqrt{b}$, for $a\ge b \ge0$ and $c\ge 0$.
As a consequence, we have that 
\[
\max_{\nu \in \{\mu,\mu'\}} \fP^{\otimes n}_{\nu_t}(\psi_{\sf LR} \neq \nu) \le  \Big(1-\frac{1}{2} H^2(\mu_t,\mu_t') \Big)^n \le  \Big(1-\frac{1}{2(1-\eta)} H^2(\tilde \mu_t, \tilde\mu_t') \Big)^n \, .
\]
By the proof of Theorem~\ref{THM:main1}, we therefore have that
\[
\max_{\nu \in \{\mu,\mu'\}} \fP^{\otimes n}_{\nu_t}(\psi_{\sf LR} \neq \nu) \le  \Big(1-\frac{1}{2(1-\eta)} H^2(\tilde \mu_t, \tilde\mu_t') \Big)^n \le \Big(1-\frac{(\eta/3)^{5/2}}{16(1-\eta)} \|\tilde\mu_t - \tilde\mu'_t\|_\pi^2 \Big)^n\, .
\]
Indeed, since $\tilde \mu$ and $\tilde \mu'$ have $\eta/3$-bounded likelihood ratios, so do $\tilde \mu_t$ and $\tilde \mu'_t$. Further, by linearity, it holds that
\[
\|\tilde\mu_t - \tilde\mu'_t\|_\pi^2 = (1-\eta)^2\|\mu_t - \mu'_t\|_\pi^2\, .
\]
As a consequence, we finally have
\[
\max_{\nu \in \{\mu,\mu'\}} \fP^{\otimes n}_{\nu_t}(\psi_{\sf LR} \neq \nu) \le \Big(1-\frac{(\eta/3)^{5/2}(1-\eta)}{16} \|\mu_t - \mu'_t\|_\pi^2 \Big)^n \le e^{-c n \|\mu_t - \mu'_t\|_\pi^2}\, ,
\]
with $c>0$ for any choice of $\eta \in (0,1)$.
\end{proof} 
\subsection{Proofs from Section~\ref{SEC:time}}

\begin{proof}[of Theorem~\ref{THM:time}] 
	Let $u_{[i]}$ be the left eigenvector of $P$ corresponding to eigenvalue $\lambda_{[i]}$. Let $\mu = \pi + \alpha u_{[i]}$ and $\mu^\prime = \pi - \alpha u_{[i]}$, where $\alpha > 0$ is sufficiently small so that $\mu$ and $\mu^\prime$ are valid probability distributions and $\mu, \mu^\prime, \pi$ pairwise satisfy the bounded likelihood ratio assumption with parameter $\varepsilon$. Let $n_0(\varepsilon, \delta) > 0$ be the sample complexity required to distinguish between $\mu$ and $\mu^\prime$ with probability greater than $\frac{1}{2} + \delta$. 

	Let $n \geq n_0$ be the sample complexity required to distinguish between distributions $\mu_t$ and $\mu^\prime_t$. Using Theorems~\ref{THM:main1} and~\ref{THM:main2}, we know that $n \asymp \lnormpi{\mu_t - \mu_t^\prime}^{-2}$ and $n_0 \asymp \lnormpi{\mu - \mu^\prime}^{-2}$. By definition, we have that $\mu - \mu^\prime = 2 \alpha u_{[i]}$ and as a result $\lnormpi{\mu - \mu^\prime} = 2 \alpha$ and $\lnormpi{\mu_t - \mu_t^\prime} = 2 \alpha \lambda_{[i]}^{t}$. Thus, we have:
	\[ \frac{n}{n_0} \asymp \left(\frac{\lnormpi{\mu - \mu^\prime}}{\lnormpi{\mu_t - \mu^\prime_t}} \right)^2 = \lambda^{-2t}_{[i]} \]
	We can invert the above to express $t$ in terms of $n$, $n_0$ and $\lambda_{[i]}$ to get the required result, which is tight up to terms involving only $\varepsilon$ and $\delta$.
\end{proof}

\subsection{Proofs from Section~\ref{sec:applications}}
\begin{proof}[Proof Sketch of Proposition~\ref{prop:random-markov}]
	We write $P = D^{-1} U$, where $D$ is a diagonal matrix with $D_{ii} =
	\sum_{x \in cX} U_{i,x}$. Writing $m_U$ to be the mean of $U_{i, j}$, which we treat to be a constant,
	standard concentration inequalities imply that $\Vert D \Vert \asymp \Vert
	D^{-1} \Vert \asymp m_U d$. It follows from \cite[Theorem
	1.2]{BCC:2010spectrum-two-examples} that $|\lambda_{[2]}(P)| \asymp
	\frac{1}{\sqrt{d}}$. 

	To establish bounds on $|\lambda_{[d]}|$, we consider $P^{-1} = U^{-1} D$
	(it is known that $P$ is non-singular with probability $1 - o(1)$), then
	$|\lambda_{[d]}|^{-1} \leq \Vert D \Vert \cdot \Vert U^{-1} \Vert$. Using a result
	of~\cite{Ver:2014invertibility}, we know that $\Vert U^{-1} \Vert =
	O(\sqrt{d})$, which gives us $|\lambda_{[d]}^{-1}| = O(m_U d^{3/2})$. 

	For a lower bound, we have $\Vert P^{-1} \Vert \geq \Vert U^{-1} \Vert \cdot
	(\Vert D^{-1} \Vert)^{-1}$. By using universality results for random
	matrices, it is known that $\Vert U^{-1} \Vert \geq \sqrt{d}$. (This can be
	achieved by subtracting a rank one matrix and using interlacing results for
	eigenvalues, cf.~\cite{Cha:2009singular}.) This together with the bounds on
	$\Vert D^{-1} \Vert$ establishes that $|\lambda_{[d]}^{-1}| = \Omega(m_U
	d^{3/2})$
\end{proof}

\begin{proof}[of Proposition~\ref{prop:regular-blockmodel}]
	In the case of two blocks, we simplify the notation a bit and assume that
	$\Delta_{1, 1} = \Delta_{2, 2} = a\cdot d$ and $\Delta_{1, 2} = \Delta_{2,
	1} = b\cdot d$ The, degree of every vertex in the graph is $(a + b)d$. Let
	$V_1 \cup V_2$ be a partition of the set of nodes into the two parts.
	Consider $v : V \rightarrow \reals$, where $v(i) = 1$ if $i \in V_1$ and
	$-1$ if $i \in V_2$. It is easily checked that $v$ is an eigenvector with
	eigenvalue $(a - b)/(a + b)$. This shows that $\lambda_{[2]} \geq \frac{a -
	b}{a + b}$.

	Now, consider a vector $v : V \rightarrow \reals$, with $v(i) =
	\frac{1}{\sqrt{2}}$, $v(j) = -\frac{1}{\sqrt{2}}$ and $v(l) = 0$ for $l \neq
	i, j$. Note that $\Vert v \Vert_2 = 1$. Let $u = Pv$, then, we have:
	\begin{align*}
		u(l) = 
		\begin{cases}
			\frac{1}{d(a + b)\sqrt{2}} & \text{if } \{l, i\} \in E \wedge \{l, j\} \not\in E \\
			-\frac{1}{d(a + b)\sqrt{2}} & \text{if } \{l, j\} \in E \wedge \{l, i\} \not\in E \\
			0 & \text{otherwise}
		\end{cases}
	\end{align*}
	If we let $N(i)$ and $N(j)$ denote the neighborhoods of $i$ and $j$ respectively, we have that,
	\begin{align*}
		\Vert u \Vert^2 &= \frac{|N(i) \Delta N(j)|}{2 d^2(a + b)^2} & |S \Delta T| \text{ denotes the symmetric difference}
	\end{align*}
	We have that $\Vert P u \Vert^2 \geq |\lambda_{[d]}|$.  If we set
	$\Delta_{1, 1} = ad = d/2 - O(1)$, there exist $i, j$ such that $|N(i)
	\Delta N(j)| = O(1)$; further, if we set $b = o(1)$, then the result
	follows.
\end{proof}

\begin{proof}[of Proposition~\ref{prop:tree}]
The matrix $P$ is symmetric, one can check directly that the given vectors are indeed eigenvectors with corresponding eigenvalues, i.e. that $P v_i = \lambda_i v_i$. Finally, we have for $2\le k \le r$,
\[
\gamma_k -\gamma_{k+1} = 2\beta_{r+1-k} - \beta_{r+2-k} > 0\, ,
\]
by order of the $\beta_\ell$. For $k=1$, we have $\gamma_1-\gamma_2 = 2\beta_r$ and $\gamma_{r+1} = \beta_0 -\beta_1 >0$.
\end{proof}

\end{document}